\documentclass[12pt]{article}
\usepackage{a4}
\usepackage{amsfonts}
\usepackage{amssymb}

\begin{document}

\title{\textbf{The motivic zeta function and its smallest poles}}
\author{Dirk Segers\thanks{Postdoctoral Fellow of the
Fund for Scientific Research - Flanders (Belgium). \newline The
authors are supported by FWO-Flanders project G.0318.06.
} \and
Lise Van Proeyen \and Willem Veys}

\date{September 17, 2012}

\maketitle

\begin{abstract}
Let $f$ be a regular function on a nonsingular complex algebraic
variety of dimension $d$. We prove a formula for the motivic zeta
function of $f$ in terms of an embedded resolution. This formula is
over the Grothendieck ring itself, and specializes to the formula of
Denef and Loeser over a certain localization. We also show that the
space of $n$-jets satisfying $f=0$ can be partitioned into locally
closed subsets which are isomorphic to a cartesian product of some
variety with an affine space of dimension $\ulcorner dn/2
\urcorner$. Finally, we look at the consequences for the poles of
the motivic zeta function.
\end{abstract}

\section{Introduction}

\noindent \textbf{(1.1)} All schemes that are considered in this
paper have base field $\mathbb{C}$.

The topological Euler-Poincar\'e characteristic $\chi$ has the
following properties on complex algebraic varieties:
$\chi(V)=\chi(V')$ if $V$ is isomorphic to $V'$, $\chi(V)=\chi(V
\setminus W) + \chi(W)$ if $W$ is a closed subset of $V$, and
$\chi(V \times W) = \chi(V)\chi(W)$. The Hodge-Deligne polynomial
of complex algebraic varieties (see (1.5))is a finer invariant
which has also these properties. The finest invariant with these
properties is the class of a variety in the Grothendieck ring.

We recall this notion. The Grothendieck ring
$K_0(\mbox{Var}_{\mathbb{C}})$ of complex algebraic varieties is the
abelian group generated by the symbols $[V]$, where $V$ is a complex
algebraic variety, subject to the relations $[V]=[V']$, if $V$ is
isomorphic to $V'$, and $[V]=[V \setminus W] + [W]$, if $W$ is
closed in $V$. One can extend the Grothendieck bracket in the
obvious way to constructible sets. The ring stucture of
$K_0(\mbox{Var}_{\mathbb{C}})$ is given by $[V] \cdot [W] := [V
\times W]$. We denote by $\mathbb{L}$ the class of the affine line,
and by $\mathcal{M}_{\mathbb{C}}$ the localization of
$K_0(\mbox{Var}_{\mathbb{C}})$ with respect to $\mathbb{L}$.

\vspace{0,5cm}

\noindent \textbf{(1.2)} Let $\mbox{Sch}$ be the category of
separated schemes of finite type over $\mathbb{C}$ and let $n \in
\mathbb{Z}_{\geq 0}$. The functor $\cdot
\times_{\mbox{\scriptsize{Spec }} \mathbb{C}} \mbox{Spec }
\mathbb{C}[t]/(t^{n+1}) : \mbox{Sch} \rightarrow \mbox{Sch}$ has a
right adjoint, which we denote by $\mathcal{L}_n$. We call
$\mathcal{L}_n(V)$ the scheme of $n$-jets of $V$ and we define an
$n$-jet on $V$ as a closed point on $\mathcal{L}_n(V)$. If $f:W
\rightarrow V$ is a morphism of schemes, then we get an induced
morphism $f_n := \mathcal{L}_n(f) : \mathcal{L}_n(W) \rightarrow
\mathcal{L}_n(V)$. For $m,n \in \mathbb{Z}_{\geq 0}$ satisfying $m
\geq n$, the canonical embeddings $V \times_{\mbox{\scriptsize{Spec
}} \mathbb{C}} \mbox{Spec } \mathbb{C}[t]/(t^{n+1}) \hookrightarrow
V \times_{\mbox{\scriptsize{Spec }} \mathbb{C}} \mbox{Spec }
\mathbb{C}[t]/(t^{m+1})$ induce canonical (projection) morphisms
$\pi_n^m : \mathcal{L}_m(V) \rightarrow \mathcal{L}_n(V)$. Because
we want $\mathcal{L}_n(V)$ to be a variety (e.g. to take its class
in the Grothendieck ring), we always endow $\mathcal{L}_n(V)$ with
its reduced structure, and we interpret the morphisms above in this
context. For more information about these constructions, see
\cite{DenefLoeserGermsofarcs} or \cite{Mustata}.

If $V$ is a closed subscheme of $\mathbb{A}^m$, the \, closed \,
points \, of \, $\mathcal{L}_n(V)$ \, are \, the $(a_{ij})_{1 \leq
i \leq m,0 \leq j \leq n} \in \mathbb{A}^{m(n+1)}$ for which
$(a_{1,0}+a_{1,1}t+\cdots+a_{1,n}t^n,
\ldots,a_{m,0}+a_{m,1}t+\cdots+a_{m,n}t^n) \in
(\mathbb{C}[t]/(t^{n+1}))^m$ satisfies the equations of $V$.
Actually, we get a set of equations of (the original possibly
nonreduced version of) $\mathcal{L}_n(V)$ by substituting an
arbitrary point of $(\mathbb{C}[t]/(t^{n+1}))^m$ in the equations
of $V$. If $V$ is an arbitrary separated scheme of finite type
over $\mathbb{C}$, we apply the construction above to the elements
of an affine cover, and we glue them together. Note that an
$n$-jet can be seen as a parameterized curve modulo $t^{n+1}$.

\vspace{0,5cm}

\noindent \textbf{(1.3)} Let $X$ be a nonsingular irreducible
algebraic variety of dimension $d$, and let $f:X \rightarrow
\mathbb{A}^1$ be a non-constant regular function. Put $V =
\mbox{div}(f)$. For each $n \in \mathbb{Z}_{\geq 0}$, we consider
the induced morphism $f_n : \mathcal{L}_n(X) \rightarrow
\mathcal{L}_n(\mathbb{A}^1)$. For $n \in \mathbb{Z}_{\geq 0}$, the
set
\[
\mathcal{X}_n := \{ \gamma \in \mathcal{L}_n(X) \mid \gamma \cdot
V=n\}
\]
is a locally closed subvariety of $\mathcal{L}_n(X)$. Note that
$\gamma \cdot V = \mbox{ord}_t(f_n(\gamma))$. The motivic zeta
function $Z(t)$ of $f:X \rightarrow \mathbb{A}^1$ is by definition
\[
Z(t) := \sum_{n \geq 0} [\mathcal{X}_n]t^n \in
\mathcal{M}_{\mathbb{C}}[[t]].
\]
In a lot of papers, there is a normalization factor
$\mathbb{L}^{-dn}$ in the $(n+1)$th term in the definition of the
motivic zeta function. Note that $Z(\mathbb{L}^{-d}t) - [X \setminus
V]$ is the naive motivic zeta function from
\cite{DenefLoesergeoarc}.

\vspace{0,5cm}

\noindent \textbf{(1.4)} We now describe a formula for $Z(t)$ in
terms of an embedded resolution. Denef and Loeser deduced it by
using motivic integration. Let $h: Y \rightarrow X$ be an embedded
resolution of $f$, i.e. $h$ is a proper birational morphism from a
nonsingular variety $Y$ such that $h$ is an isomorphism on $Y
\setminus h^{-1}(f^{-1}\{0\})$ and $h^{-1}(f^{-1}\{0\})$ is a
normal crossings divisor. Let $E_i$, $i \in S$, be the irreducible
components of $h^{-1}(f^{-1}\{0\})$. Let $K_{Y/X}$ be the relative
canonical divisor supported in the exceptional locus of $h$. We
define the numerical data $N_i$ and $\nu_i$ by the equalities
$\mbox{div}(f \circ h) = \sum_{i \in S} N_i E_i$ and $K_{Y/X} =
\sum_{i \in S} (\nu_i-1)E_i$. For $I \subset S$, denote
$E_I^{\circ} := (\cap_{i \in I} E_i) \setminus (\cup_{i \notin I}
E_i)$. Then, the announced formula for $Z(t)$ is
\[
Z(t) = \sum_{I \subset S} [E_I^{\circ}] \prod_{i \in I}
\frac{(\mathbb{L}-1)
\mathbb{L}^{dN_i-\nu_i}t^{N_i}}{1-\mathbb{L}^{dN_i-\nu_i}t^{N_i}}.
\]
In particular, $Z(t)$ is rational and belongs to the subring of
$\mathcal{M}_{\mathbb{C}}[[t]]$ generated by
$\mathcal{M}_{\mathbb{C}}$ and the elements
$t^{N}/(1-\mathbb{L}^{dN-\nu}t^{N})$, with $\nu,N \in
\mathbb{Z}_{>0}$. Note that the equation implies, in particular,
that the right hand side does not depend on the embedded resolution.

\vspace{0,5cm}

\noindent \textbf{(1.5)} We now introduce the Hodge zeta function.
Recall that the Hodge-Deligne polynomial of a complex algebraic
variety $W$ is
\[
H(W) := \sum_{p,q} \left( \sum_{i \geq 0} (-1)^i
h^{p,q}\left(H_c^i(W,\mathbb{C})\right) \right)u^pv^q \in
\mathbb{Z}[u,v],
\]
where $h^{p,q}\left(H_c^i(W,\mathbb{C})\right)$ is the dimension
of the $(p,q)$-Hodge component of the $i$th cohomology group with
compact support of $W$. The Hodge zeta function of $f$ is
\[
Z_{\mathrm{Hod}}(t) := \sum_{I \subset S} H(E_I^{\circ}) \prod_{i
\in I} \frac{(uv-1) (uv)^{dN_i - \nu_i} t^{N_i}}{1- (uv)^{dN_i -
\nu_i} t^{N_i}}.
\]
The right hand side does not depend on the embedded resolution
because it is obtained from the right hand side of the formula in
(1.4) by specializing to Hodge-Deligne polynomials. Note that the
Hodge zeta function is in a lot of papers a normalization of this
one.

\vspace{0,5cm}

\noindent \textbf{(1.6)} Let $f$ and $V$ be as in (1.3). We consider
the power series
\[
J(t) := \sum_{n \geq 0} [\mathcal{L}_n(V)] t^n \in
\mathcal{M}_{\mathbb{C}}[[t]].
\]
Because $[\mathcal{X}_n] = \mathbb{L}^d [\mathcal{L}_{n-1}(V)] -
[\mathcal{L}_n(V)]$ for $n \geq 1$ and $[\mathcal{X}_0] = [X] -
[V]$, we have the relation
\[
J(t) = \frac{Z(t)-[X]}{\mathbb{L}^dt-1}.
\]
Consequently, the series $J(t)$ and $Z(t)$ determine each other.

\vspace{0,5cm}

\noindent \textbf{(1.7)} In Section 2, we prove the formula of
(1.4) without using motivic integration. We will actually prove a
stronger result which is over $K_0(\mbox{Var}_{\mathbb{C}})$
instead of $\mathcal{M}_{\mathbb{C}}$: for an integer $c$
satisfying $(\nu_i-1)/N_i \leq c$ for all $i \in S$, we have
\[
\sum_{n \geq 0} [\mathcal{X}_n] (\mathbb{L}^{2cd+c-d} t)^n =
\sum_{I \subset S} [E_I^{\circ}] \prod_{i \in I}
\frac{(\mathbb{L}-1)
\mathbb{L}^{(2cd+c)N_i-\nu_i}t^{N_i}}{1-\mathbb{L}^{(2cd+c)N_i-\nu_i}t^{N_i}}
\]
in $K_0(\mbox{Var}_{\mathbb{C}})[[t]]$. After localizing with
respect to $\mathbb{L}$, we can indeed deduce the formula of (1.4)
because $\mathbb{L}$ is not a zero-divisor in
$\mathcal{M}_{\mathbb{C}}$. However, it is unknown whether
$\mathbb{L}$ is a zero-divisor in $K_0(\mbox{Var}_{\mathbb{C}})$.
This implies that we cannot deduce our formula straightforward from
the one in (1.4) and that we do not know whether the formula of
(1.4) holds in $K_0(\mbox{Var}_{\mathbb{C}})[[t]]$ whenever $dN_i -
\nu_i \geq 0$ for all $i \in S$. (There always exists an embedded
resolution for which this condition is satisfied.)

In Section 3, we will prove that $[\mathcal{L}_n(V)]$ is a multiple
of $\mathbb{L}^{\ulcorner dn/2 \urcorner}$ in
$K_0(\mbox{Var}_{\mathbb{C}})$ for all $n \in \mathbb{Z}_{\geq 0}$.
(We use the notation $\ulcorner x \urcorner$ for the smallest
integer larger than or equal to $x \in \mathbb{R}$.) We will
actually construct a partition of $\mathcal{L}_n(V)$ into locally
closed subsets which are isomorphic to $W \times
\mathbb{A}^{\ulcorner dn/2 \urcorner}$ for some variety $W$
depending on the locally closed subset. The first author proved
already an analogous result for the number of solutions of
polynomial congruences in \cite{Segersmathann}. The difficulty here
is that we do not have to count solutions, but that we have to
construct isomorphisms. We also note that our setting of (1.3) is
more general than polynomials, i.e. regular functions on affine
space.

In Section 4, we will consider $Z(t)$ as a power series over a ring
$R$ which is a quotient of the image of the localization map in
$\mathcal{M}_{\mathbb{C}}$. Using the previous result, we prove that
$Z(t)$ belongs to the subring of $R[[t]]$ generated by $R[t]$ and
the elements $1/(1-\mathbb{L}^{dN-\nu}t^{N})$, with $\nu,N \in
\mathbb{Z}_{>0}$ and $\nu/N \leq d/2$. An analogous result was
already proved in \cite{Segersmathann} for the topological zeta
function (and for Igusa's $p$-adic zeta function), where it says
that there are no poles (with real part) less than $-d/2$. See
\cite{RodriguesVeys} for a possible definition of the notion of pole
for the motivic zeta function. Because the ring $R$ specializes to
Hodge-Deligne polynomials, this result is also true for the Hodge
zeta function.

In Section 5, we adapt the previous results in a relative setting.

\vspace{0,5cm}

\noindent \textsl{Acknowledgements.} We want to thank the referee
for the careful reading of the paper and for the useful remarks and
suggestions.

\section{The motivic zeta function over the Grothendieck ring}

\noindent \textbf{(2.1)} If $g : Y \rightarrow Z$ is \'etale, one
can sometimes reduce a problem about $\mathcal{L}_m(Y)$ to an
analogous problem for $\mathcal{L}_m(Z)$ because of the following
proposition. For a proof, see for example \cite[Proposition
2.2]{Blickle}.

\vspace{0,2cm}

\noindent \textbf{Proposition.} \textsl{Let $g : Y \rightarrow Z$
be \'etale and $m \in \mathbb{Z}_{>0}$. Then the natural map
$\mathcal{L}_m(Y) \rightarrow Y \times_Z \mathcal{L}_m(Z)$ is an
isomorphism.}

\vspace{0,2cm}

If $g : Y \rightarrow \mathbb{A}^d$ is \'etale, we obtain
\[
\mathcal{L}_m(Y) \cong Y \times_{\mathbb{A}^d}
\mathcal{L}_m(\mathbb{A}^d) = Y \times_{\mathbb{A}^d} ( \mathbb{A}^d
\times \mathbb{A}^{dm}) \cong Y \times \mathbb{A}^{dm}.
\]
If $Y$ is an arbitrary nonsingular irreducible algebraic variety of
dimension $d$, we can cover $Y$ with open subsets $U$ for which
$\mathcal{L}_m(U) \cong U \times \mathbb{A}^{dm}$. Consequently,
$[\mathcal{L}_m(Y)] = [Y] \mathbb{L}^{dm}$.

Note that also the equality $[\mathcal{X}_n] = \mathbb{L}^d
[\mathcal{L}_{n-1}(V)] - [\mathcal{L}_n(V)]$ for $n \geq 1$ of (1.6)
can be proved by using this proposition.

\vspace{0,5cm}

\noindent \textbf{(2.2)} We also need a theorem of Denef and
Loeser \cite[Lemma 3.4]{DenefLoeserGermsofarcs} to obtain our
formula.

\vspace{0,2cm}

\noindent \textbf{Theorem.} \textsl{Let $X$ and $Y$ be nonsingular
irreducible algebraic varieties of dimension $d$, and let $h: Y
\rightarrow X$ be a proper birational morphism. For $e,m \in
\mathbb{Z}_{\geq 0}$ satisfying $m \geq e$, the set
\[
\Delta_{e,m} := \{ \gamma \in \mathcal{L}_m(Y) \mid \gamma \cdot
K_{Y/X}=e\}
\]
is a locally closed subset of $\mathcal{L}_m(Y)$. If $m \geq 2e$,
then $\Delta_{e,m}$ is the union of fibers of $h_m$ and the
restriction $\Delta_{e,m} \rightarrow h_m(\Delta_{e,m})$ of $h_m$ is
a piecewise trivial fibration with fiber $\mathbb{A}^e$. Moreover,
two elements of the same fiber have the same image in
$\mathcal{L}_{m-e}(Y)$.}

\vspace{0,5cm}

\noindent \textbf{(2.3)} Let $X$ be a nonsingular irreducible
algebraic variety of dimension $d$, and let $f:X \rightarrow
\mathbb{A}^1$ be a non-constant regular function. Let $h: Y
\rightarrow X$ be an embedded resolution of $f$. Let $E_i$, $i \in
S$, be the irreducible components of $h^{-1}(f^{-1}\{0\})$, and
let $N_i$ and $\nu_i$, with $i \in S$, be the numerical data.

\vspace{0,2cm}

\noindent \textbf{Theorem.} \textsl{If $c$ is an integer
satisfying $(\nu_i-1)/N_i \leq c$ for all $i \in S$, then
\begin{eqnarray} \label{formulegrothendieckring}
\sum_{n \geq 0} [\mathcal{X}_n]
(\mathbb{L}^{2cd+c-d} t)^n = \sum_{I \subset S} [E_I^{\circ}]
\prod_{i \in I} \frac{(\mathbb{L}-1)
\mathbb{L}^{(2cd+c)N_i-\nu_i}t^{N_i}}{1-\mathbb{L}^{(2cd+c)N_i-\nu_i}t^{N_i}}
\end{eqnarray}
in $K_0(\mbox{Var}_{\mathbb{C}})[[t]]$.}

\vspace{0,2cm}

\noindent \textsl{Proof.} For $n,m \in \mathbb{Z}_{\geq 0}$
satisfying $m \geq n$, the set
\[
\mathcal{X}_{n,m} := \{ \gamma \in \mathcal{L}_m(X) \mid \gamma
\cdot \mbox{div}(f)=n\}
\]
is a locally closed subvariety of $\mathcal{L}_m(X)$. Note that
$\mathcal{X}_n = \mathcal{X}_{n,n}$. Because $\mathcal{X}_{n,m}
\cong \mathcal{X}_n \times \mathbb{A}^{d(m-n)}$ if $X$ admits an
\'etale map $X \rightarrow \mathbb{A}^d$, we have
$[\mathcal{X}_{n,m}] = [\mathcal{X}_n] \mathbb{L}^{d(m-n)}$ for
general $X$.

Let $\gamma \in h_m^{-1}(\mathcal{X}_{n,m})$. We have that
\[
\sum_{i \in S} N_i (\gamma \cdot E_i) = \gamma \cdot (\sum_{i \in
S} N_i E_i) = h_m(\gamma) \cdot \mbox{div}(f) = n,
\]
and that
\[
\gamma \cdot K_{Y/X} = \gamma \cdot (\sum_{i \in S} (\nu_i-1) E_i)
= \sum_{i \in S} (\nu_i-1) (\gamma \cdot E_i).
\]
Let $c$ be an integer satisfying $(\nu_i-1)/N_i \leq c$ for all $i
\in S$. Such an integer exists because $S$ is finite. Note that
there exists an embedded resolution $h$ (which is a composition of
blowing-ups with well chosen centra) for which $\nu_i/N_i \leq
d-1$ for all $i \in S$ \cite[Proof of Theorem
2.4.0]{Segersthesis}. We obtain
\[
\gamma \cdot K_{Y/X} = \sum_{i \in S} (\nu_i-1) (\gamma \cdot E_i)
\leq c \sum_{i \in S} N_i (\gamma \cdot E_i) = cn.
\]
In particular, there are only a finite number of possibilities for
$\gamma \cdot K_{Y/X}$. In view of (2.2), we will try to find a
formula for
\[
\sum_{n \geq 0} \mathbb{L}^{cn} [\mathcal{X}_{n,2cn}]t^n = \sum_{n
\geq 0} [\mathcal{X}_n] (\mathbb{L}^{2cd+c-d}t)^n \in
K_0(\mbox{Var}_{\mathbb{C}})[[t]]
\]
in terms of $h$.

Fix $n \in \mathbb{Z}_{>0}$. For an $|S|$-tuple of positive
integers $a=(a_i)_{i \in S}$, we define $S_a:=\{i \in S \mid a_i
> 0\}$, $E_a^{\circ} := E_{S_a}^{\circ}$ and $|a| :=
\mbox{Card}(S_a)$. We define a set $A$ by
\[
A := \{ (a_i)_{i \in S} \mid \forall i \in S \, : \, a_i \in
\mathbb{Z}_{\geq 0} \; , \; \sum_{i \in S} a_i N_i = n \mbox{ and
} E_a^{\circ} \not= \emptyset \} ,
\]
and obtain a disjoint union
\[
h^{-1}_{2cn}(\mathcal{X}_{n,2cn}) = \bigsqcup_{a \in A} \{ \gamma
\in \mathcal{L}_{2cn}(Y) \mid \forall i \in S \, : \, \gamma \cdot
E_i = a_i\}.
\]

Fix $a \in A$. Put $e = \sum_{i \in S} a_i (\nu_i-1)$. Denote the
origin of $\gamma$ by $\gamma_0$. If $U_t$, $t \in T$, is a
partition of $E_a^{\circ}$, then
\[
\{ \gamma \in \mathcal{L}_{2cn}(Y) \mid \forall i \in S \, : \,
\gamma \cdot E_i = a_i\} = \bigsqcup_{t \in T} \{ \gamma \in
\mathcal{L}_{2cn}(Y) \mid \forall i \in S \, : \, \gamma \cdot E_i
= a_i \mbox{ and } \gamma_0 \in U_t \}.
\]
Using (2.1), we can take \cite[Proof of Proposition 2.5]{Craw} a
partition $U_t$, $t \in T$, of $E_a^{\circ}$ into locally closed
subset such that for each $t \in T$, the set
\[
F_{a,t} := \{ \gamma \in \mathcal{L}_{2cn}(Y) \mid \forall i \in S
\, : \, \gamma \cdot E_i = a_i \mbox{ and } \gamma_0 \in U_t \}
\]
is isomorphic to $U_t \times \mathbb{A}^{2cdn - \sum_{i \in S_a}
a_i} \times (\mathbb{A}^1 \setminus \{0\})^{|a|}$, and hence
\[
[F_{a,t}] = [U_t] \mathbb{L}^{2cdn - \sum_{i \in S_a} a_i}
(\mathbb{L}-1)^{|a|}.
\]
Because $F_{a,t} \subset \Delta_{e,2cn}$ is a union of fibres of
$h_{2cn}$ (one can use the last statement in Theorem 2.2 to prove
this), we obtain from (2.2) that the restriction $F_{a,t}
\rightarrow h_{2cn}(F_{a,t})$ of $h_{2cn}$ is a piecewise trivial
fibration with fibre $\mathbb{A}^e$. Hence,
\[
[F_{a,t}] = \mathbb{L}^e [h_{2cn}(F_{a,t})].
\]
By summing over all $t \in T$, we obtain
\begin{eqnarray*}
\mathbb{L}^{e} [h_{2cn}(\{ \gamma \in \mathcal{L}_{2cn}(Y) \mid
\forall i \in S \, : \, \gamma \cdot E_i = a_i \})] & = &
\mathbb{L}^e [h_{2cn}(\bigsqcup_{t \in T} F_{a,t})] \\ & = &
\sum_{t \in T} \mathbb{L}^e [h_{2cn}(F_{a,t})] \\ & = & \sum_{t \in T} [F_{a,t}] \\
& = & [E_a^{\circ}] \mathbb{L}^{2cdn - \sum_{i \in S_a} a_i}
(\mathbb{L}-1)^{|a|}.
\end{eqnarray*}

Now we multiply both sides of the obtained equality with
$\mathbb{L}^{cn-e}$ and sum over all $a \in A$. Note that $e$
depends on $a$. We obtain
\begin{eqnarray*}
\mathbb{L}^{cn} [\mathcal{X}_{n,2cn}] & = & \sum_{a \in A}
[E_a^{\circ}] \mathbb{L}^{2cdn-\sum_{i \in S_a} a_i}
\mathbb{L}^{cn-e} (\mathbb{L}-1)^{|a|} \\ & = & \sum_{a \in A}
[E_a^{\circ}]
(\mathbb{L}-1)^{|a|} \mathbb{L}^{(2cd+c)n-e-\sum_{i \in S_a} a_i} \\
& = & \sum_{a \in A} [E_a^{\circ}] (\mathbb{L}-1)^{|a|}
\mathbb{L}^{\sum_{i \in S_a} ((2cd+c)N_i - \nu_i)a_i} \\ & = &
\sum_{a \in A} [E_a^{\circ}] \prod_{i \in S_a} (\mathbb{L}-1)
\mathbb{L}^{((2cd+c)N_i - \nu_i)a_i}.
\end{eqnarray*}
The last expression is the coefficient of $t^n$ in the formal
power series
\[
\sum_{I \subset S} [E_I^{\circ}] \prod_{i \in I}
\frac{(\mathbb{L}-1)
\mathbb{L}^{(2cd+c)N_i-\nu_i}t^{N_i}}{1-\mathbb{L}^{(2cd+c)N_i-\nu_i}t^{N_i}},
\]
and consequently, we have finished our proof. $\qquad \Box$

\vspace{0,2cm}

\noindent \textsl{Remark.} If we consider
(\ref{formulegrothendieckring}) as an equality in
$\mathcal{M}_{\mathbb{C}}[[t]]$ and if we replace $t$ by
$\mathbb{L}^{-(2cd+c-d)}t$, we obtain
\[
Z(t) = \sum_{I \subset S} [E_I^{\circ}] \prod_{i \in I}
\frac{(\mathbb{L}-1)
\mathbb{L}^{dN_i-\nu_i}t^{N_i}}{1-\mathbb{L}^{dN_i-\nu_i}t^{N_i}}
\]
in $\mathcal{M}_{\mathbb{C}}[[t]]$. This is the formula of (1.4)
which was first proved by Denef and Loeser using motivic
integration.

\section{Divisibility of jet spaces in $K_0(\mbox{Var}_{\mathbb{C}})$}

\noindent \textbf{(3.1)} Let $X$ be a $d$-dimensional complex
analytic manifold with analytic coordinates $(u_1,\ldots,u_d)$.
These coordinates induce tangent vector fields $\partial/\partial
u_1,\ldots,$ $\partial/\partial u_d$ along $X$. Let $f$ be a complex
analytic function on $X$ and let $b \in X$. By Taylor's theorem, we
have for points $x$ in a small enough neighbourhood of $b$ that
\begin{eqnarray*}
f(x) & = & \lim_{m \rightarrow \infty} \sum_{\alpha \in
\mathbb{Z}_{\geq 0}^d : |\alpha| \leq m}
\frac{(\partial^{|\alpha|}f/
\partial u^{\alpha})(b)}{\alpha!} (x-b)^{\alpha}   \\
& = & f(b) + \sum_{j=1}^d (\partial f/\partial u_j)(b) (x_j-b_j) +
\cdots
\end{eqnarray*}
We explain the notation. The coordinates of $b$ are
$(b_1,\ldots,b_d)$ and those of $x$ are $(x_1,\dots,x_d)$. For
$\alpha=(\alpha_1,\ldots,\alpha_d) \in \mathbb{Z}_{\geq 0}^d$, we
put $\alpha! = \alpha_1! \cdots \alpha_d!$, $(x-b)^{\alpha} =
(x_1-b_1)^{\alpha_1} \cdots (x_d-b_d)^{\alpha_d}$, $|\alpha| =
\alpha_1 + \cdots + \alpha_d$ and $\partial^{|\alpha|}f/\partial
u^{\alpha} =
\partial^{\alpha_1+\cdots+\alpha_d}f/\partial u_1^{\alpha_1} \cdots \partial
u_d^{\alpha_d}$.

Let $\gamma$ be a convergent arc in $X$, i.e. a $d$-tuple of
convergent power series in $t$ with origin in $X$. Let $l \in
\mathbb{Z}_{>0}$. Write $\gamma = b+t^lz$, where $b$ is a $d$-tuple
of polynomials in $t$ of degree less than $l$ and where $z$ is a
$d$-tuple of convergent power series. For every $t \in \mathbb{C}$
in the convergence domain of $\gamma$ (and for which $\gamma(t) \in
X$), we have
\[
f(\gamma(t)) = \lim_{m \rightarrow \infty} \sum_{\alpha \in
\mathbb{Z}_{\geq 0}^d : |\alpha| \leq m}
\frac{(\partial^{|\alpha|}f/
\partial u^{\alpha})(b(t))}{\alpha!} t^{|\alpha|l} (z(t))^{\alpha},
\]
and consequently, we obtain the equality
\[
f(\gamma) = \lim_{m \rightarrow \infty} \sum_{\alpha \in
\mathbb{Z}_{\geq 0}^d : |\alpha| \leq m}
\frac{(\partial^{|\alpha|}f/
\partial u^{\alpha})(b)}{\alpha!} t^{|\alpha|l} z^{\alpha}
\]
as formal power series in $t$.

Since every $n$-jet is liftable to a convergent arc, we get for an
$n$-jet of the form $\gamma=b+t^lz$ that
\[
f(\gamma) = \sum_{|\alpha|=0}^{\ulcorner n/l \urcorner}
\frac{(\partial^{|\alpha|}f/
\partial u^{\alpha})(b)}{\alpha!} t^{|\alpha|l} z^{\alpha} \mbox{ mod
}t^{n+1}.
\]

\vspace{0,5cm}

\noindent \textbf{(3.2)} Let $X$ be a nonsingular irreducible
algebraic variety of dimension $d$ and let $g:X \rightarrow
\mathbb{A}^d$ be an \'etale map. We identify $\mathcal{L}_n(X)$ and
$X \times_{\mathbb{A}^d} \mathcal{L}_n(\mathbb{A}^d)$ by using the
canonical isomorphism of (2.1).

The coordinates $(u_1,\ldots,u_d)$ on $\mathbb{A}^d$ induce
analytic coordinates on a complex neighbourhood of every point of
$X$. The tangent vector fields $\partial/\partial
u_1,\ldots,\partial/\partial u_d$, which we define along the whole
of $X$, are actually algebraic. This implies that all first and
higher order partial derivatives of a regular function on $X$ with
respect to $u_1,\ldots,u_d$ are regular functions on $X$.

Let $f$ be a regular function on $X$. Let $l \in \{1,\ldots,n\}$.
Let $(x,b) \in X \times_{\mathbb{A}^d} \mathcal{L}_n(\mathbb{A}^d) =
\mathcal{L}_n(X)$ and let $z=(z_1,\ldots,z_d) \in
\mathcal{L}_{n-l}(\mathbb{A}^d)$. Then
\[
f(x,b+t^l z) = \sum_{|\alpha|=0}^{\ulcorner n/l \urcorner}
\frac{(\partial^{|\alpha|}f/
\partial u^{\alpha})(x,b)}{\alpha!} t^{|\alpha|l} z^{\alpha} \mbox{ mod } t^{n+1}.
\]

\vspace{0,5cm}

\noindent \textbf{(3.3) Theorem.} \textsl{Let $X$ be a nonsingular
irreducible algebraic variety of dimension $d \in \mathbb{Z}_{>1}$
and let $f:X \rightarrow \mathbb{A}^1$ be a regular function. Put
$V := \mbox{div}(f)$. Then $[\mathcal{L}_n(V)]$ is a multiple of
$\mathbb{L}^{\ulcorner dn/2 \urcorner}$ in
$K_0($Var$_{\mathbb{C}})$ for all $n \in \mathbb{Z}_{\geq 0}$.}

\vspace{0,2cm}

\noindent \textsl{Remark.} It follows that $[\mathcal{X}_n]$ is also
a multiple of $\mathbb{L}^{\ulcorner dn/2 \urcorner}$ for all $n \in
\mathbb{Z}_{\geq 0}$.

\vspace{0,2cm}

\noindent \textsl{Proof.} We are going to partition
$\mathcal{L}_n(V)$ into a finite number of locally closed subsets
which are isomorphic to $W \times \mathbb{A}^{\ulcorner dn/2
\urcorner}$ for some variety $W$ depending on the locally closed
subset.

If the theorem holds for the members of an open cover of $X$ and all
their intersections, then it holds for $X$, by additivity of the
Grothendieck bracket. Hence, we may assume that there exists an
\'etale map $g:X \rightarrow \mathbb{A}^d$.

Let $r$ be $n/2$ if $n$ is even and $(n+1)/2$ if $n$ is odd. We
define
\[
\mathcal{L}_{n,r}(V) = \left\{ (x,b) \in \mathcal{L}_n(X) \left|
\begin{array}{l} \forall j \in \{1,\ldots,d\} \, : \, (\partial f / \partial
u_j)(x,b) \equiv 0 \mbox{ mod } t^r \mbox{ and} \\
f(x,b) \equiv 0 \mbox{ mod } t^{n+1}
\end{array} \right. \right\}
\]
and for $k \in \{0,1,\ldots,r-1\}$, we define
\[
\mathcal{L}_{n,k}(V) = \left\{ (x,b) \in \mathcal{L}_n(X) \left|
\begin{array}{l} \forall j \in \{1,\ldots,d\} \, : \, (\partial f /
\partial u_j)(x,b) \equiv 0 \mbox{ mod } t^k, \\
\exists j \in \{1,\ldots,d\} \, : \, (\partial f / \partial
u_j)(x,b) \not\equiv 0 \mbox{ mod } t^{k+1} \mbox{ and}
\\ f(x,b) \equiv 0 \mbox{ mod } t^{n+1}
\end{array} \right. \right\}.
\]
Then, the sets $\mathcal{L}_{n,k}(V)$, $k \in \{0,1,\ldots,r\}$,
are locally closed subsets of $\mathcal{L}_n(V)$ which partition
$\mathcal{L}_n(V)$.

\vspace{0,3cm}

We prove that $\mathcal{L}_{n,r}(V)$ is isomorphic to $W \times
\mathbb{A}^{\ulcorner dn/2 \urcorner}$ for some variety $W$. Let
$(x,b) \in \mathcal{L}_{n,r}(V)$. For $z=(z_1,\ldots,z_d) \in
\mathcal{L}_{r-1}(\mathbb{A}^d)$, we have that
\begin{eqnarray*}
f(x,b+t^{n-r+1}z) & = & f(x,b) + \sum_{j=1}^d (\partial f /
\partial u_j)(x,b) t^{n-r+1} z_j + t^{2(n-r+1)} (
\ldots ) \\ & \equiv & 0 \mbox{ mod } t^{n+1},
\end{eqnarray*}
such that $(x,b + t^{n-r+1} \mathcal{L}_{r-1}(\mathbb{A}^d)) \subset
\mathcal{L}_n(V)$. Because $n-r+1 \geq r$, we obtain that $(x,b +
t^{n-r+1} \mathcal{L}_{r-1}(\mathbb{A}^d)) \subset
\mathcal{L}_{n,r}(V)$. Consequently, $\mathcal{L}_{n,r}(V) \cong
\pi^n_{n-r}(\mathcal{L}_{n,r}(V)) \times \mathbb{A}^{dr}$. This
proves our assertion for $\mathcal{L}_{n,r}(V)$ because $dr \geq
\ulcorner dn/2 \urcorner$.

\vspace{0,3cm}

Let $k \in \{0,1,\ldots,r-1\}$. We study $\mathcal{L}_{n,k}(V)$.
Let $p \in \{ 1,\ldots,d\}$ and let $l \in \{k,\ldots,n-k\}$. We
define
\[
\mathcal{L}_{n,k,p}(V) = \left\{ (x,b) \in \mathcal{L}_n(X) \left|
\begin{array}{l} \forall j \in \{1,\ldots,d\} \, : \, (\partial f /
\partial u_j)(x,b) \equiv 0 \mbox{ mod } t^k \\ \mbox{and }
(\partial f / \partial u_p)(x,b) \not\equiv 0 \mbox{ mod } t^{k+1}
\\ \mbox{and } f(x,b) \equiv 0 \mbox{ mod } t^{n+1}
\end{array} \right. \right\}.
\]

\noindent and

\[
\mathcal{O}_{l,k,p}(V) = \left\{ (x,b) \in \mathcal{L}_l(X) \left|
\begin{array}{l} \forall j \in \{1,\ldots,d\} \, : \, (\partial f /
\partial u_j)(x,b) \equiv 0 \mbox{ mod } t^k \\ \mbox{and }
(\partial f / \partial u_p)(x,b) \not\equiv 0 \mbox{ mod } t^{k+1}
\\ \mbox{and } f(x,b) \equiv 0 \mbox{ mod } t^{k+l+1}
\end{array} \right. \right\}.
\]

\noindent Note that $\mathcal{L}_{n,0,p}(V) =
\mathcal{O}_{n,0,p}(V)$. We check that $f(x,b)$ is well defined
modulo $t^{k+l+1}$ in the definition of $\mathcal{O}_{l,k,p}(V)$.
Suppose that $(x,b) \in \mathcal{L}_l(X)$ satisfies $(\partial f /
\partial u_j)(x,b) \equiv 0 \mbox{ mod } t^k$ for all
$j \in \{1,\ldots,d\}$. Then
\begin{eqnarray*}
f(x,b+t^{l+1}z) & = & f(x,b) + \sum_{j=1}^d (\partial f /
\partial u_j)(x,b) t^{l+1} z_j + t^{2(l+1)} ( \ldots )
\\ & \equiv & f(x,b) \mbox{ mod } t^{k+l+1},
\end{eqnarray*}
and consequently, $f(x,b)$ is well defined modulo $t^{k+l+1}$.

\noindent The following isomorphisms can be checked easily:
\begin{eqnarray*}
\mathcal{L}_{n,k,p}(V) & \cong & \mathcal{O}_{n-k,k,p}(V) \times
\mathbb{A}^{dk}, \\ \mathcal{O}_{l+1,k,p}(V) & \cong &
\mathcal{O}_{l,k,p}(V) \times \mathbb{A}^{d-1}, \mbox{ for } l \in
\{k,\ldots,n-k-1\},
\end{eqnarray*}
where the projections on the first factor are respectively
$\pi^n_{n-k}$ and $\pi^{l+1}_l$. The isomorphism
$\mathcal{O}_{l,k,p}(V) \times \mathbb{A}^{d-1} \rightarrow
\mathcal{O}_{l+1,k,p}(V)$ maps
$((x;b_{10}+b_{11}t+\cdots+b_{1l}t^l,\ldots,b_{d0}+b_{d1}t+\cdots+b_{dl}t^l),
(a_1,\ldots,a_{d-1}))$ to the unique element of
$\mathcal{O}_{l+1,k,p}(V)$ of the form
$(x;b_{10}+b_{11}t+\cdots+b_{1l}t^l+a_1t^{l+1},\ldots,
b_{p0}+b_{p1}t+\cdots+b_{pl}t^l+ct^{l+1},\ldots,
b_{d0}+b_{d1}t+\cdots+b_{dl}t^l+a_{d-1}t^{l+1})$. The unique value
$c \in \mathbb{C}$ is obtained by solving a linear equation.

\noindent Consequently,
\[
\mathcal{L}_{n,k,p}(V) \cong \mathcal{O}_{k,k,p}(V) \times
\mathbb{A}^{(d-1)(n-2k)+dk},
\]
where the projection on the first factor is $\pi^n_k$. This
implies
\begin{eqnarray*}
\mathcal{L}_{n,k,1}(V) & \cong & \mathcal{O}_{k,k,1}(V) \times
\mathbb{A}^{(d-1)(n-2k)+dk}, \\
\mathcal{L}_{n,k,2}(V) \setminus \mathcal{L}_{n,k,1}(V) & \cong &
(\mathcal{O}_{k,k,2}(V) \setminus \mathcal{O}_{k,k,1}(V)) \times
\mathbb{A}^{(d-1)(n-2k)+dk}, \\ & \vdots & \\
\mathcal{L}_{n,k,d}(V) \setminus (\cup_{1\leq p \leq d-1}
\mathcal{L}_{n,k,p}(V)) & \cong & \left( \mathcal{O}_{k,k,d}(V)
\setminus (\cup_{1\leq p \leq d-1} \mathcal{O}_{k,k,p}(V)) \right)
\times \mathbb{A}^{(d-1)(n-2k)+dk}.
\end{eqnarray*}
This finishes our proof because the left hand sides form a
partition of $\mathcal{L}_{n,k}(V)$ into locally closed subsets
and because $(d-1)(n-2k)+dk \geq \ulcorner dn/2 \urcorner$.
$\qquad \Box$

\section{The smallest poles of motivic zeta functions}

\noindent \textbf{(4.1)} Let $X$ be a nonsingular irreducible
algebraic variety of dimension $d \in \mathbb{Z}_{>1}$ and let
$f:X \rightarrow \mathbb{A}^1$ be a regular function. Put $V :=
\mbox{div}(f)$. In this section, we fix an embedded resolution for
which $dN_i-\nu_i \geq 0$ for every $i \in S$. Note that we
mentioned already in (2.3) that there always exists an embedded
resolution which satisfies the stronger condition $\nu_i/N_i \leq
d-1$ for every $i \in S$.

Consider the ideal
\[
I' = \{ \alpha \in K_0(\mbox{Var}_{\mathbb{C}}) \mid \exists k \in
\mathbb{Z}_{\geq 0} \, : \, \mathbb{L}^k \alpha = 0\}
\]
of $K_0(\mbox{Var}_{\mathbb{C}})$ and put $R' =
K_0(\mbox{Var}_{\mathbb{C}})/I'$. The image of $\mathbb{L}$ in $R'$,
which is also denoted by $\mathbb{L}$, is not a zero-divisor in
$R'$. Note that $I'$ is the kernel of the localization map
$K_0(\mbox{Var}_{\mathbb{C}}) \rightarrow \mathcal{M}_{\mathbb{C}}$,
such that $R'$ can be considered as the image of this map in
$\mathcal{M}_{\mathbb{C}}$. Consequently, the formula of (1.4) still
holds if we consider $Z(t)$ as a power series over $R'$.

We want to use for instance that $\cap_{k \in \mathbb{Z}_{\geq 0}}
(\mathbb{L}^k) = \{ 0 \}$ and that a number $k \in \mathbb{Z}
\setminus \{0\}$ is not a zero-divisor. Because we do not know
whether these are true in $R'$, we will work in an appropriate
quotient of $R'$. Consider the ideal
\begin{eqnarray*}
I & = & \bigcap_{k \in \mathbb{Z}_{\geq 0}} \{ \alpha \in R' \mid
\exists n \in \mathbb{Z} \setminus \{0\} \, : \, n \alpha \in
(\mathbb{L}^k) \} \\ & = & \bigcap_{k \in \mathbb{Z}_{\geq 0}} \{
\alpha \in R' \mid \exists n \in \mathbb{Z} \setminus \{0\} \, : \,
n \alpha \mbox{ is divisible by } \mathbb{L}^k \mbox{ in } R' \}
\end{eqnarray*}
of $R'$ and put $R = R'/I$. Note that $R$ specializes to
Hodge-Deligne polynomials and that we do not know whether $I \not=
\{0\}$. One verifies easily that an element of $\mathbb{Z} \setminus
\{0\}$ is not a zero-divisor in $R$. One also checks that the image
in $R$ of $\{ \alpha \in R' \mid \exists n \in \mathbb{Z} \setminus
\{0\} \, : \, n \alpha \in (\mathbb{L}^k) \}$ contains
$(\mathbb{L}^k)$, and consequently $\cap_{k \in \mathbb{Z}_{\geq 0}}
(\mathbb{L}^k) = \{ 0 \}$ in $R$. Thus, if $\alpha$ is a non-zero
element of $R$, there exists a $k \in \mathbb{Z}_{\geq 0}$ such that
$\alpha$ is divisible by $\mathbb{L}^k$ but not by
$\mathbb{L}^{k+1}$ in $R$. Moreover, if $\alpha$ is a non-zero
element of $R$, there exists a positive integer $k$ which has for
every $n \in \mathbb{Z} \setminus \{0\}$ the property that $n
\alpha$ is not divisible by $\mathbb{L}^k$. We also have that
$1-\mathbb{L}^k$, with $k \in \mathbb{Z}_{>0}$, is not a
zero-divisor in $R$. Indeed, if $\alpha \in R$ satisfies
$(1-\mathbb{L}^k)\alpha=0$, then $\alpha = \mathbb{L}^k \alpha =
\mathbb{L}^{2k} \alpha = \mathbb{L}^{3k} \alpha = \cdots$, and thus
$\alpha \in \cap_{k \in \mathbb{Z}_{\geq 0}} (\mathbb{L}^k) = \{ 0
\}$.

From now on, we will consider the motivic zeta function $Z(t)$ as a
power series over $R$. The formula of $Z(t)$ in terms of an embedded
resolution also holds over $R$. We write the motivic zeta function
in the form
\[
Z(t) = \frac{B(t)}{\prod_{i \in I}
(1-\mathbb{L}^{dN_i-\nu_i}t^{N_i})},
\]
where $I \subset S$ and where $B(t)$ is not divisible by any of the
$1-\mathbb{L}^{dN_i-\nu_i}t^{N_i}$, with $i \in I$. Put $l:=\min
\{-\nu_i/N_i \mid i \in I \}$.

In the next paragraphs, we will work in a more general context. By
abuse of notation, we will use the symbols of this particular
situation.

\vspace{0,5cm}

\noindent \textbf{(4.2)} Let $Z(t)$ be an \textsl{arbitrary} element
of $R[[t]]$ of the form
\[
Z(t) = \frac{B(t)}{\prod_{i \in I}
(1-\mathbb{L}^{dN_i-\nu_i}t^{N_i})},
\]
where every $(\nu_i,N_i) \in \mathbb{Z}_{>0} \times
\mathbb{Z}_{>0}$ satisfies $dN_i-\nu_i \geq 0$ and where $B(t) \in
R[t]$ is not divisible by any of the
$1-\mathbb{L}^{dN_i-\nu_i}t^{N_i}$, with $i \in I$. Put $l:=\min
\{-\nu_i/N_i \mid i \in I \}$. Define the elements $\gamma_n \in
R$ by the equality
\[
Z(t) = \sum_{n \geq 0} \gamma_n t^n.
\]

\vspace{0,5cm}

\noindent \textbf{(4.3) Proposition.} \textsl{There exists an
integer $a$ which is independent of $n$ such that $\gamma_n$ is a
multiple of $\mathbb{L}^{\ulcorner(d+l)n-a\urcorner}$ in $R$ for all
integers $n$ satisfying $(d+l)n-a \geq 0$.}

\vspace{0,2cm}

\noindent \textsl{Remark.} (i) The statement in the proposition is
obviously equivalent to the following. If $l' \leq l$, then there
exists an integer $a$ which is independent of $n$ such that
$\gamma_n$ is a multiple of
$\mathbb{L}^{\ulcorner(d+l')n-a\urcorner}$ for all integers $n$
satisfying $(d+l')n-a \geq 0$.
\\(ii) Suppose that we are in the situation of (4.1). It follows
from (4.6) that $d+l>0$, so that $(d+l)n-a$ rises linearly as a
function of $n$ with a slope depending on $l$. The condition
$(d+l)n-a \geq 0$ is thus satisfied for $n$ large enough.

\vspace{0,2cm}

\noindent \textsl{Proof.} We will say that a formal power series in
$t$ has the divisibility property if the coefficient of $t^n$ is a
multiple of $\mathbb{L}^{\ulcorner(d+l)n\urcorner}$ for every $n$.

For $i \in I$, the series
\[
\frac{1}{1-\mathbb{L}^{dN_i-\nu_i}t^{N_i}}= \sum_{n \geq 0}
\mathbb{L}^{n(dN_i-\nu_i)} t^{nN_i}
\]
has the divisibility property because $dN_i-\nu_i$ is an integer
larger than or equal to $N_i(d+l)$.

One can easily check that the product of a finite number of power
series with the divisibility property also has the divisibility
property. Let $g$ be the degree of $B(t)$. For $n \geq g$, we will
have that $\gamma_n$ is a multiple of $\mathbb{L}^{\ulcorner
(d+l)(n-g) \urcorner}$. This implies our statement. $\qquad \Box$

\vspace{0,5cm}

\noindent \textbf{(4.4)} We will decompose $Z(t)$ into partial
fractions in (4.6). To this end, we need to apply the following
lemma several times.

\vspace{0,2cm}

\noindent\textbf{Lemma.} \textsl{(a) Let $i,j \in I$ such that
$\nu_i/N_i \not= \nu_j/N_j$. Then, there exist polynomials
$g(x,t),h(x,t) \in \mathbb{Z}[x,t]$ and an integer $k \in
\mathbb{Z}_{>0}$ such that
\[
g(x,t) (1-x^{dN_i-\nu_i}t^{N_i}) + h(x,t) (1-x^{dN_j-\nu_j}t^{N_j})
= 1-x^k
\]
holds in $\mathbb{Z}[x,t]$, and consequently such that
\[
g(\mathbb{L},t) (1-\mathbb{L}^{dN_i-\nu_i}t^{N_i}) + h(\mathbb{L},t)
(1-\mathbb{L}^{dN_j-\nu_j}t^{N_j}) = 1-\mathbb{L}^k
\]
holds in $R[t]$. \\ (b) Let $D(t) \in R[t]$. There exist
polynomials $g(t), h(t) \in R[t]$ with $\mbox{deg}(h)<N_i$ and a
$k \in \mathbb{Z}_{\geq 0}$ such that}
\[
\mathbb{L}^k D(t) = (1-\mathbb{L}^{dN_i-\nu_i}t^{N_i}) g(t) + h(t).
\]

\vspace{0,2cm}

\noindent \textsl{Proof.} (a) Although $\mathbb{Z}[x,t]$ is not a
PID, we can obtain the first relation by applying the algorithm of
Bezout-Bachet in number theory to the polynomials
$1-x^{dN_i-\nu_i}t^{N_i}$ and $1-x^{dN_j-\nu_j}t^{N_j}$ in the
variable $t$. The number $k$ is different from 0 because otherwise
the polynomials $1-x^{dN_i-\nu_i}t^{N_i}$ and
$1-x^{dN_j-\nu_j}t^{N_j}$ would have a non-trivial common divisor,
and this is not the case because $\nu_i/N_i \not= \nu_j/N_j$. \\ (b)
This is straightforward by applying the division algorithm.  $\qquad
\Box$

\vspace{0,5cm}

\noindent \textbf{(4.5)} For $r \in \mathbb{Z}_{>0}$, we define a
function $f_r : \mathbb{Z}_{\geq 0} \rightarrow \mathbb{Z}_{\geq 0}$
by the relation
\[
\frac{1}{(1-x)^r} = \sum_{n=0}^{\infty} f_r(n) x^n.
\]
One proves by induction on $r$ that
\[
f_r(n) = \frac{(n+r-1)!}{n!(r-1)!} =
\frac{(n+1)(n+2)\ldots(n+r-1)}{(r-1)!}.
\]

\vspace{0,2cm}

\noindent \textbf{Lemma.} \textsl{Let $m \in \mathbb{Z}_{>1}$. Let
$n_1,\ldots,n_m$ be $m$ different natural numbers. Then, the
determinant of the matrix with the elements
\begin{eqnarray*}
v_1 & = & (f_1(n_1),f_2(n_1),\ldots,f_m(n_1)) \\  v_2 & = &
(f_1(n_2),f_2(n_2),\ldots,f_m(n_2)) \\ & \vdots & \\ v_m & = &
(f_1(n_m),f_2(n_m),\ldots,f_m(n_m))
\end{eqnarray*}
of $\mathbb{Z}^m$ in the rows is equal to}
\[
\frac{\prod_{j>i} (n_j-n_i)}{\prod_{i=1}^{m-1} i!}.
\]

\vspace{0,2cm}

\noindent \textsl{Remark.} (i) This determinant is different from
zero, so the set $\{v_1,v_2,\ldots,v_m\}$ is linearly independent.
Consequently, every element $e_i$ of the standard basis of the
$\mathbb{Z}$-module $\mathbb{Z}^m$ has a multiple which is generated
by it. \\ (ii) This lemma is probably known. We include its proof by
lack of reference.

\vspace{0,2cm}

\noindent \textsl{Proof.} The proof is by induction on $m$. The
statement is trivial for $m=2$. Let now $m>2$. We expand the
determinant along the last column, apply the induction hypothesis to
the cofactors, use Vandermonde determinants, and obtain that it is
equal to
\[ \frac{1}{\prod_{i=1}^{m-2} i!}
\left| \begin{array}{cccccc}
1 & n_1 & n_1^2 & \cdots & n_1^{m-2} & f_m(n_1) \\
1 & n_2 & n_2^2 & \cdots & n_2^{m-2} & f_m(n_2) \\
\vdots & \vdots & \vdots & & \vdots & \vdots \\
1 & n_m & n_m^2 & \cdots & n_m^{m-2} & f_m(n_m)
\end{array} \right|. \]
By using properties of determinants and the Vandermonde determinant,
we see that this is equal to
\[ \frac{1}{\prod_{i=1}^{m-1} i!}
\left| \begin{array}{cccccc}
1 & n_1 & n_1^2 & \cdots & n_1^{m-2} & n_1^{m-1} \\
1 & n_2 & n_2^2 & \cdots & n_2^{m-2} & n_2^{m-1} \\
\vdots & \vdots & \vdots & & \vdots & \vdots \\
1 & n_m & n_m^2 & \cdots & n_m^{m-2} & n_m^{m-1}
\end{array} \right| = \frac{\prod_{j>i} (n_j-n_i)}{\prod_{i=1}^{m-1} i!}. \qquad \Box \]

\vspace{0,5cm}

\noindent \textbf{(4.6) Proposition.} \textsl{There exist an integer
$a$ which is independent of $n$ and positive integers $N$ and $b$
such that $\gamma_{nN+b}$ is not a multiple of
$\mathbb{L}^{\ulcorner (d+l)(nN+b)+a \urcorner}$ in $R$ for $n$
large enough.}

\vspace{0,2cm}

\noindent \textsl{Proof.}  Put $I_1=\{j \in I \mid -\nu_i/N_i=l\}$
and $I_2=I \setminus I_1$. Let $N$ be the lowest common multiple of
the $N_i$, $i \in I_1$, and let $\nu$ be the lowest common multiple
of the $\nu_i$, $i \in I_1$. Remark that $\nu/N=\nu_i/N_i$ for all
$i \in I_1$. Let $m$ be the cardinality of $I_1$. Because
$1-\mathbb{L}^{dN-\nu}t^N$ is a multiple of
$1-\mathbb{L}^{dN_i-\nu_i}t^{N_i}$ for all $i \in I_1$, we can write
\[
Z(t) = \frac{D(t)}{(1-\mathbb{L}^{dN-\nu}t^N)^m \prod_{i \in I_2}
(1-\mathbb{L}^{dN_i-\nu_i}t^{N_i})},
\]
where $D(t) \in R[t]$. Applying decomposition into partial fractions
(see Lemma 4.4), we can write
\begin{eqnarray}
wZ(t) & = &
\frac{\mu_{m,0}+\mu_{m,1}t+\cdots+\mu_{m,N-1}t^{N-1}}{(1-\mathbb{L}^{dN-\nu}t^N)^m}
+\frac{\mu_{m-1,0}+\mu_{m-1,1}t+\cdots+\mu_{m-1,N-1}t^{N-1}}{(1-\mathbb{L}^{dN-\nu}t^N)^{m-1}}
\nonumber \\ & & + \cdots +
\frac{\mu_{1,0}+\mu_{1,1}t+\cdots+\mu_{1,N-1}t^{N-1}}{1-\mathbb{L}^{dN-\nu}t^N}
+ \frac{E(t)}{\prod_{i \in I_2} (1-\mathbb{L}^{dN_i-\nu_i}t^{N_i})}
\nonumber \\ & = & \sum_{b=0}^{N-1} \sum_{n=0}^{\infty} (f_m(n)
\mu_{m,b} + \cdots + f_1(n) \mu_{1,b}) \mathbb{L}^{ndN-n\nu}
t^{nN+b} \label{deel1}
\\ & & + \frac{E(t)}{\prod_{i \in I_2}
(1-\mathbb{L}^{dN_i-\nu_i}t^{N_i})}, \label{deel2}
\end{eqnarray}
where $\mu_{i,j} \in R$, where $E(t) \in R[t]$ and where $w$ is a
product of elements of the form $1-\mathbb{L}^k$ and $\mathbb{L}^k$,
with $k > 0$. Note that $wD(t)$ is not divisible by
$1-\mathbb{L}^{dN-\nu}t^N$ because $w$ is not a zero divisor in $R$,
the constant term of $1-\mathbb{L}^{dN-\nu}t^N$ is a unit in $R$ and
$D(t)$ is not divisible by $1-\mathbb{L}^{dN-\nu}t^N$.

We now consider the first part (\ref{deel1}) of $wZ(t)$. Because
$wD(t)$ is not divisible by $(1-\mathbb{L}^{dN-\nu}t^N)^m$, there
exists a $b \in \{0,\ldots,N-1\}$ for which the coefficient of
$t^{nN+b}$ is different from 0 for infinitely many $n$. Fix from
now on such a $b$ and a $j \in \{1,\ldots,m\}$ for which
$\mu_{j,b} \not= 0$. Take a positive integer $c$ such that we have
for every $n \in \mathbb{Z} \setminus \{0\}$ that $n \mu_{j,b}$ is
not divisible by $\mathbb{L}^c$. There do not exist $m$ positive
integers $n_1,\ldots,n_m$ for which $f_m(n_1) \mu_{m,b} + \cdots +
f_1(n_1) \mu_{1,b},\ldots,f_m(n_m) \mu_{m,b} + \cdots + f_1(n_m)
\mu_{1,b}$ are multiples of $\mathbb{L}^c$, because otherwise, we
can use Lemma 4.5 to obtain that $\mu_{j,b}$ has an integer
multiple which is a multiple of $\mathbb{L}^c$. Consequently, for
$n$ large enough, $f_m(n) \mu_{m,b} + \cdots + f_1(n) \mu_{1,b}$
is not a multiple of $\mathbb{L}^c$. The coefficient of $t^{nN+b}$
in the power series expansion of (\ref{deel1}) is equal to
$(f_m(n) \mu_{m,b} + \cdots + f_1(n) \mu_{1,b})
\mathbb{L}^{(d+l)nN}$, which is not a multiple of
$\mathbb{L}^{(d+l)nN+c} = \mathbb{L}^{(d+l)(nN+b)-(d+l)b+c}$ for
$n$ large enough. So let $a$ be the largest integer smaller than
or equal to $c-(d+l)b$.

Now we consider the remaining part (\ref{deel2}) of $wZ(t)$. We
obtain from Proposition 4.3 that there exists an $l'>l$ and an
integer $a'$ such that the coefficient of $t^n$ in the power
series expansion of (\ref{deel2}) is a multiple of
$\mathbb{L}^{\ulcorner(d+l')n-a'\urcorner}$ for $n$ large enough.
Consequently, this coefficient is a multiple of
$\mathbb{L}^{\ulcorner(d+l)n+a\urcorner}$ for $n$ large enough.

Because $w\gamma_{nN+b}$ is the sum of two elements of which exactly
one is a multiple of $\mathbb{L}^{\ulcorner(d+l)(nN+b)+a\urcorner}$
for $n$ large enough, we obtain that $w\gamma_{nN+b}$, and thus also
$\gamma_{nN+b}$, is not a multiple of
$\mathbb{L}^{\ulcorner(d+l)(nN+b)+a\urcorner}$ for $n$ large enough.
$\qquad \Box$

\vspace{0,2cm}

\noindent \textsl{Corollaries.} (i) If there exists an integer $a$
such that $\gamma_n$ is a multiple of $\mathbb{L}^{\ulcorner
(d+l')n-a \urcorner}$ for all $n$ satisfying $(d+l')n-a \geq 0$,
then $l' \leq l$. This is the converse of Proposition 4.3. \\ (ii)
Because we saw in the previous section that $[\mathcal{X}_n]$ is a
multiple of $\mathbb{L}^{\ulcorner dn/2 \urcorner}$ if we are in the
situation of (4.1), we obtain that $l \geq -d/2$.

\vspace{0,2cm}

Because of the second corollary, we have proved the following
theorem.

\vspace{0,2cm}

\noindent \textbf{Theorem.} \textsl{The motivic zeta function
$Z(t) \in R[[t]]$ belongs to the subring of $R[[t]]$ generated by
$R[t]$ and the elements $1/(1-\mathbb{L}^{dN-\nu}t^{N})$, with
$\nu,N \in \mathbb{Z}_{>0}$ and $\nu/N \leq d/2$.}

\vspace{0,5cm}

\noindent \textbf{(4.7)} In (4.1), we denoted the image of the
localization map $K_0(\mbox{Var}_{\mathbb{C}}) \rightarrow
\mathcal{M}_{\mathbb{C}}$ by $R'$. We introduced an ideal $I$ of
$R'$ and put $R=R'/I$. The previous theorem is a priori weaker
than the analogous statement over $R'$ (or
$\mathcal{M}_{\mathbb{C}}$), but it is not if $I=\{0\}$. We do not
know whether $I \not= \{0\}$, but at any rate the theorem
specializes to Hodge-Deligne polynomials. This gives us the
following.

\vspace{0,2cm}

\noindent \textbf{Theorem.} \textsl{The Hodge zeta function
$Z_{\mathrm{Hod}}(t)$ belongs to the subring of $\mathbb{Q}(u,v)(t)$
generated by $\mathbb{Q}(u,v)[t]$ and the elements $1/(1-(uv)^{dN -
\nu}t^N)$, with $\nu,N \in \mathbb{Z}_{>0}$ and $\nu/N \leq d/2$.}

\section{The relative setting}

The generalization to the relative setting was suggested by the
referee. Let $X$ be a nonsingular irreducible algebraic variety of
dimension $d$, and let $f: X \rightarrow \mathbb{A}^1$ be a
non-constant regular function. Let $X_0$ be the reduced scheme
determined by $f=0$. Note that $X_0 = \mathcal{L}_0(V)$, where
$V=\mbox{div}(f)$ as before. For $n \geq 1$, we have that
$\mathcal{X}_n$ is an $X_0$ variety because of the canonical
morphism $\pi^n_0: \mathcal{X}_n \rightarrow X_0$. Therefore, we can
consider the class $[\mathcal{X}_n / X_0]$ of $\mathcal{X}_n$ in the
relative Grothendieck ring $K_0(\mbox{Var}_{X_0})$ of
$X_0$-varieties. The definition is the straightforward
generalization of the usual one, see for example
\cite{DenefLoesergeoarc}.

One obtains analogously as in section 2 that
\begin{eqnarray*}
\sum_{n \geq 1} [\mathcal{X}_n / X_0] (\mathbb{L}^{2cd+c-d} t)^n =
\sum_{\emptyset \not= I \subset S} [E_I^{\circ} / X_0] \prod_{i
\in I} \frac{(\mathbb{L}-1)
\mathbb{L}^{(2cd+c)N_i-\nu_i}t^{N_i}}{1-\mathbb{L}^{(2cd+c)N_i-\nu_i}t^{N_i}}
\end{eqnarray*}
in $K_0(\mbox{Var}_{X_0})[[t]]$. Here, $c$ is an arbitrary integer
satisfying $(\nu_i-1)/N_i \leq c$ for all $i \in S$ and
$\mathbb{L}$ is the class of $\mathbb{A}^1 \times X_0$ in
$K_0(\mbox{Var}_{X_0})$.

The main result of section 3 can also be adapted to this context.
Suppose that the dimension $d$ of $X$ is in $\mathbb{Z}_{>1}$. Then
$[\mathcal{L}_n(V) / X_0]$ is a multiple of $\mathbb{L}^{\ulcorner
dn/2 \urcorner}$ in $K_0(\mbox{Var}_{X_0})$ for all $n \in
\mathbb{Z}_{\geq 0}$.

Also section 4 can be generalized. One analogously constructs a
ring $R$ from $K_0(\mbox{Var}_{X_0})$ such that $Z(t)$, considered
as an element of $R[[t]]$, belongs to the subring of $R[[t]]$
generated by $R[t]$ and the elements
$1/(1-\mathbb{L}^{dN-\nu}t^{N})$, with $\nu,N \in \mathbb{Z}_{>0}$
and $\nu/N \leq d/2$.

\footnotesize{

\noindent \textsl{Address:} University of Leuven \\
\mbox{}\hspace{1,37cm} Department of Mathematics
\\ \mbox{}\hspace{1,37cm} Celestijnenlaan 200B \\
\mbox{}\hspace{1,37cm} B-3001 Leuven \\ \mbox{}\hspace{1,37cm} Belgium

\vspace{0,2cm}

\noindent\textsl{E-mail:} dirk.segers@wis.kuleuven.be \\
\mbox{}\hspace{1,16cm} lise.vanproeyen@wis.kuleuven.be \\
\mbox{}\hspace{1,16cm} wim.veys@wis.kuleuven.be

\vspace{0,2cm}

\noindent \textsl{URL:} http://wis.kuleuven.be/algebra/segers/segers.htm} \\
\mbox{} \hspace{0,73cm} http://wis.kuleuven.be/algebra/veys.htm

\end{document}